\documentclass[11pt]{amsart}
\usepackage[paper=a4paper,left=30mm,width=150mm,top=25mm,bottom=25mm]{geometry}
 \usepackage{amsthm}
 \usepackage{comment}
 \usepackage{amssymb}
 \usepackage{amsmath}
 \usepackage{amsfonts}
 \usepackage{mathrsfs}
 \usepackage{amscd}
 \usepackage[all,cmtip]{xy}
 \author{Kevin Ventullo}
 \title{On the rank one abelian Gross-Stark conjecture}
\date{}
\newtheorem{conjecture}{Conjecture}
\newtheorem{thm}{Theorem}
\newtheorem{lma}{Lemma}

\newtheorem{prp}{Proposition}

\newcommand{\Z}{\mathbb{Z}}
\newcommand{\Zp}{\mathbb{Z}_p}
\newcommand{\Q}{\mathbb{Q}}
\newcommand{\Qp}{\mathbb{Q}_p}
\newcommand{\C}{\mathbb{C}}

\begin{document}

  \begin{abstract}
  Let $F$ be a totally real number field, $p$ a rational prime, and $\chi$ a finite order totally odd abelian character of Gal$(\overline{F}/F)$ such that $\chi(\mathfrak{p})=1$ for some $\mathfrak{p}|p$. Motivated by a conjecture of Stark, Gross conjectured a relation between the derivative of the $p$-adic $L$-function associated to $\chi$ at its exceptional zero and the $\mathfrak{p}$-adic logarithm of a $p$-unit in the $\chi$ component of $F_\chi^\times$. In a recent work, Dasgupta, Darmon, and Pollack have proven this conjecture assuming two conditions: that Leopoldt's conjecture holds for $F$ and $p$, and that if there is only one prime of $F$ lying above $p$, a certain relation holds between the $\mathscr{L}$-invariants of $\chi$ and $\chi^{-1}$. The main result of this paper removes both of these conditions, thus giving an unconditional proof of the conjecture.
  \end{abstract}

\maketitle
\tableofcontents
  \section{Introduction}

   Let $F$ be a totally real field of degree $g>1$. Fix a prime $p$ and embeddings $\overline{\Qp}\hookleftarrow \overline{\Q}\hookrightarrow \C$. Let $\chi:\text{Gal}(\overline{F}/F)\rightarrow \overline{\Q}^\times$ be a totally odd character of conductor $\mathfrak{n}$, and $F_\chi$ the cyclic extension of $F$ cut out by $\chi$. Let $\omega:\text{Gal}(\Q(\mu_{2p})/\Q)\rightarrow (\Z/2p\Z)^\times$ denote the Teichm$\ddot{\text{u}}$ller character. Let $S$ be any finite set of primes of $F$ including all archimedean primes. Associated to $\chi$ and $S$ is a complex analytic function $L_S(\chi,s)$ defined for Re$(s)>1$ by

   \begin{equation*}
   L_S(\chi,s)=\sum_{(\mathfrak{a},S)=1}\chi(\mathfrak{a})N(\mathfrak{a})^{-s}=\prod_{\mathfrak{p}\not\in S}(1-\chi(\mathfrak{p})N\mathfrak{p}^{-s})^{-1},
   \end{equation*}
   which has a holomorphic continuation to all of $\mathbb{C}$. By Siegel's rationality theorem, $L_S(\chi,1-k)\in\overline{\Q}$ for $k\geq 1$. Using the functional equation, one can show that the order of vanishing of $L_{S}(\chi,s)$ at $s=0$ is equal to the number of $v\in S$ such that $\chi(v)=1$.

   Let us now assume the set $S$ contains all places above $p$. Let $F_\infty$ be the cyclotomic $\Zp$-extension of $F$, and $\Gamma=\text{Gal}(F_\infty/F)$, which is canonically isomorphic to a subgroup of $1+\hookrightarrow 1+2p\Zp$. For use later, we fix a topological generator $u$ of $\Gamma$, which gives an isomorphism $\Zp[[\Gamma]]\cong \Zp[[T]]=:\Lambda$ via $u\mapsto 1+T$. We will identify $u$ with its image in $1+2p\Zp$. A character of $\text{Gal}(\overline{F}/F)$ is said to be of type $S$, resp. type $W$, if the extension it cuts out is disjoint from $F_\infty$, resp. contained in $F_\infty$. Since $\Gamma$ is a direct summand of Gal$(F^{ab}/F)$, any character can be decomposed as a product of a type $S$ character and a type $W$ character, which we write as $\chi=\chi_S\chi_W$.

   By work of Deligne and Ribet~\cite{DR:80}, there is a pseudo-measure $\mathcal{L}_{S,\chi\omega}\in\text{Frac}(\Zp[[\Gamma]])$ which interpolates classical $L$-values via the formula
   \begin{equation*}
   \chi_{cyc}^k\omega^{-k}\psi(\mathcal{L}_{S,\chi\omega}) = L_{S}(\psi\chi\omega^{1-k},1-k).
   \end{equation*}
    where $\psi$ is any character of type $W$. We also use $\mathcal{L}_{S,\chi\omega}$ to denote the corresponding element of $F_\Lambda:=\text{Frac}(\Lambda)$ via the isomorphism above. Then the previous formula can be written
    \begin{equation*}
    \mathcal{L}_{S,\chi\omega}(\zeta u^k -1) = L_{S}(\psi\chi\omega^{1-k},1-k),
    \end{equation*}
    where $\zeta=\psi(u)$. Taking $\zeta=1$, we get a $p$-adic analytic function
    \begin{align*}
    L_{p,S}(\chi\omega):\Zp&\rightarrow\overline{\Qp}\\s&\mapsto \mathcal{L}_{S,\chi\omega}(u^{1-s}-1).
    \end{align*}
    This is usually referred to as the $p$-adic $L$-function of the even character $\chi\omega$. We see from the above properties that it satisfies
    \begin{equation*}
    L_{p,S}(\chi\omega,-n)=L_{S}(\chi\omega^{-n},-n)\text{ for all }n\geq 0.
    \end{equation*}

    It follows that if $L_{S}(\chi,s)$ vanishes at $s=0$, then so does $L_{p,S}(\chi\omega,s)$. In this setting, Gross has formulated a conjecture which one can think of as comparing the $p$-adic derivative of the left hand side with the archimedean derivative of the right hand side. Suppose that $ord_{s=0}L_S(\chi,s)=1$. Then there is a unique $\mathfrak{p}\in S$ such that $\chi(\mathfrak{p})=1$. If $\mathfrak{p}\nmid p$, there is a simple relation between the respective derivatives coming from the relation between $L_{p,S\backslash\{\mathfrak{p}\}}(\chi\omega,s)$ and $L_{S\backslash\{p\}}(\chi\omega,s)$. If $\mathfrak{p}|p$, a more sophisticated construction is required.

    Let $E$ be a finite extension of $\Qp$ containing all values of $\chi$. Let $\mathfrak{P}$ be a prime of $F_\chi$ lying over $\mathfrak{p}$. By our assumptions on $\chi$ and $S$, the subspace $U_\chi$ of $\mathcal{O}_{F_\chi,S}^\times\otimes E$ on which Gal$(F_\chi/F)$ acts by $\chi^{-1}$ is one-dimensional over $E$. Let $0\neq u_\chi\in U_\chi$. Define
    \begin{equation*}
    \mathscr{L}_{alg}(\chi):=\frac{((\log_p\circ Norm_{F_{\chi,\mathfrak{P}}/\Qp})\otimes id) (u_\chi)}{(ord_\mathfrak{P}\otimes id)(u_\chi)}\in E.
    \end{equation*}

    In ~\cite{Gr}, Gross conjectures, and proves for $F=\mathbb{Q}$, the following
    \begin{conjecture}
    Let $F$ be a totally real number field, $p$ a prime, $\chi$ a finite order character of $F$ such that $\chi(\mathfrak{p})=1$ for some $\mathfrak{p}|p$, $S$ the set of primes of $F$ dividing $\text{cond}(\chi)p\infty$, and $R=S\setminus \{\mathfrak{p}\}$. Then\\
    \indent i.\phantom{i} If $L_R(\chi,0)=0$, then $L_{p,S}'(\chi\omega,0)=0$.\\
    \indent ii. If $L_R(\chi,0)\neq 0$, then $L_{p,S}'(\chi\omega,0)=\mathscr{L}_{alg}(\chi)L_R(\chi,0)$.
    \end{conjecture}

    The first part of this conjecture follows from the stronger statement that the order of vanishing of the $p$-adic $L$-function at an exceptional zero is greater than or equal to that of the archimedean $L$-function\footnote{In fact, Gross conjectures that these orders of vanishing are equal.}. In ~\cite[Lemma 1.2]{DDP:09}, this stronger statement is shown to follow from the Iwasawa Main Conjecture for the character $\chi$ (they assume $\chi$ is of type $S$, but this is actually not necessary; see Lemma 1 below). Unfortunately, the proof of the Main Conjecture is not quite complete when $p=2$ (see ~\cite[\S 11]{Wi:90}). However, the inequality between orders of vanishing has recently been shown for all $p$ by Speiss ~\cite{Sp} and Charollois-Dasgupta ~\cite{CD:12}, by entirely different methods. Thus, Conjecture $1.i$ is known in all cases.

In [DDP], Conjecture 1.\emph{ii} is proven under the following assumptions:
\begin{itemize}
\item Leopoldt's Conjecture is true for $F$ and $p$.
\item If $\mathfrak{p}$ is the unique prime above $p$, then
\begin{equation*}
    ord_{k=1}(\mathscr{L}_{an}(\chi,k)+\mathscr{L}_{an}(\chi^{-1},k))=ord_{k=1}\mathscr{L}_{an}(\chi^{-1},k),
\end{equation*}
\end{itemize}
where
\begin{align*}
\mathscr{L}_{an}(\chi,k):= \frac{-L_{p,S}(\chi\omega,1-k)}{L_R(\chi,0)}\\
\mathscr{L}_{an}(\chi):=\frac{d}{dk}\mathscr{L}_{an}(\chi,k)|_{k=1}.
\end{align*}

In this paper we prove
\begin{thm}
Conjecture 1 is true unconditionally.
\end{thm}

In the above notation, Gross's conjecture can be stated as $\mathscr{L}_{an}(\chi)=\mathscr{L}_{alg}(\chi)$. Let $d_\chi$ denote the order of vanishing of $\mathscr{L}_{an}(\chi,k)$ at $k=1$, and similarly for $d_{\chi^{-1}}$. Then Dasgupta-Darmon-Pollack's second condition is equivalent to assuming $d_\chi\geq d_{\chi^{-1}}$, and if they are equal, the leading terms of $\mathscr{L}_{an}(\chi,k)$ and $\mathscr{L}_{an}(\chi^{-1},k)$ at $k=1$ shouldn't cancel.

To remove Leopoldt's conjecture, we construct in Section 4 a certain ordinary family of parallel weight Hilbert modular forms with weight zero specialization equal to the constant form $\textbf{1}$ (see Theorem 2). In Section 2, we recall the proof of Theorem 1 given in \emph{loc. cit.}, but assuming the existence of this family so as to remove Leopoldt from the hypotheses. In Section 3, we remove the condition on $\mathscr{L}$-invariants by breaking into two cases: first assuming $d_\chi<d_{\chi^{-1}}$, and second assuming $d_\chi=d_{\chi^{-1}}$ and the leading terms cancel. At the end of Section 4, we use Theorem 2 to give a simplified proof of the ``Leopoldt'' part of the Iwasawa Main Conjecture.\\
\\
\indent \emph{Acknowledgements}. I am grateful to Samit Dasgupta, Henri Darmon, and Rob Pollack for their beautifully written paper to which this paper owes its existence. I am especially grateful to Samit for several helpful conversations, encouragement, and for suggesting the method by which we construct the $\Lambda$-adic form in Section 4.

I am also grateful to Haruzo Hida for answering my questions and providing helpful comments.

Finally, I want to thank Chandrashekhar Khare for his guidance and support, as well as suggesting a careful reading of Dasgupta-Darmon-Pollack's work.\\
\\

\section{Dasgupta-Darmon-Pollack's Proof}

\subsection{Conjecture 1.i}
We begin by showing the Iwasawa Main Conjecture implies part $i$. of Conjecture 1.
\begin{lma}
Let $\chi$ be a finite order character of $F$. If IMC holds for $(\chi_S,p)$ (e.g. if $p>2$) then
\begin{equation*}
   ord_{s=0}L_{p,S}(\chi\omega,s)\geq ord_{s=0}L_{S}(\chi,s).
\end{equation*}
   \begin{proof}
   We may assume $S$ is minimal, i.e. consists only of the primes dividing $cond(\chi)p\infty$. Let $d=ord_{s=0}L_{S}(\chi,s)$, which, by minimality of $S$, is just the number of primes $\mathfrak{p}|p$ in $F$ for which $\chi(\mathfrak{p})=1$. Let $\zeta=\chi_W(u)$. Then
   \begin{equation*}
   ord_{s=0}L_{p,S}(\chi\omega,s)=ord_{T=u-1}\mathcal{L}_{\chi\omega}=ord_{T=\zeta u-1}\mathcal{L}_{\chi_S\omega}.
   \end{equation*}

Let $F_{\chi,\infty^-}$ be the maximal anticyclotomic $\Zp$-extension of $F_\chi$, and $F_{\chi,\infty^-,ur}$ the maximal subextension which becomes unramified over $F_{\chi,\infty}$. The following fact seems to be well known, but we include a proof here for completeness:
\begin{equation*}
  dim_E (\text{Gal}(F_{\chi,\infty^-,ur}/F_{\chi})\otimes_{\Zp} E)^{\chi^{-1}} = d.
\end{equation*}

     To see this, note that by class field theory,
     \begin{equation*}
     (\text{Gal}(F_{\chi,\infty^-}/F_{\chi})\otimes_{\Zp} E)^{\chi^{-1}}\cong \prod_{\mathcal{O}_F\supset\mathfrak{p}|p} \left(\prod _{\mathcal{O}_{F_\chi}\supset \mathfrak{P}|\mathfrak{p}} U_\mathfrak{P}\otimes_{\Zp} E\right)^{\chi^{-1}},
     \end{equation*}

    where $U_\mathfrak{P}$ denotes the units in the $\mathfrak{P}$-adic completion of $\mathcal{O}_{F_\chi}$, and the superscript denotes the $\chi^{-1}$ component of this space as a $\text{Gal}(F_\chi/F)$-module. The maximal quotient of this group which becomes unramified over $F_{\chi,\infty}$ is obtained by taking the quotient of each $U_\mathfrak{P}$ by the kernel of the norm map to $\Zp^\times$. We can write this as
    \begin{equation*}
    \prod_{\mathcal{O}_F\supset\mathfrak{p}|p} \left(\prod _{\mathcal{O}_{F_\chi}\supset \mathfrak{P}|\mathfrak{p}} E\right)^{\chi^{-1}},
     \end{equation*}
     where $\text{Gal}(F_{\chi}/F)$ acts by permuting the factors inside the parentheses for each $\mathfrak{p}$. Since $\chi^{-1}$ is by definition a faithful character of $\text{Gal}(F_{\chi}/F)$, it will appear as a constituent of the expression inside the parentheses precisely when $\mathfrak{p}$ splits completely in $F_\chi/F$, i.e. $\chi(\mathfrak{p})=1$.

   To complete the proof of the lemma, note that the subgroup Gal$(F_{\chi_S,\infty}/F_\infty)\cong\text{Gal}(F_{\chi_S}/F)$ acts on this extension by $\chi_S^{-1}$, and $\text{Gal}(F_{\chi_S,\infty}/F_{\chi_S})\cong\text{Gal}(F_\infty/F)$ acts by $\chi_W^{-1}$. It follows from the Main Conjecture that $d\leq ord_{T=\zeta u-1}\mathcal{L}_{\chi_S\omega}$. This implies the lemma.
   \end{proof}
   \end{lma}

\subsection{Classical and $\Lambda$-adic Hilbert Modular Forms}
Fix $F,p, \mathfrak{n}$ as above. Let $\mathfrak{d}$ be the different of $F$. Let $U_F$ denote the units of $\mathcal{O}_F$, and $U_F^+$ the totally positive units. Let $\mathfrak{c}$ be a representative of a strict ideal class in $\mathcal{O}_F$, and $\mathfrak{c}^+$ the cone of positive elements. Let $\varphi:(\mathcal{O_F}/\mathfrak{n})^\times \rightarrow \overline{\Q}^\times$ be a character.

\textbf{Definition.} A complex $\mathfrak{c}$-Hilbert modular form of weight $k$, level $\mathfrak{n}$, and character $\psi$ is a holomorphic function $f$ on the product of $g$ upper half planes, indexed by the embeddings of $F$ into $\mathbb{R}$, such that for every element of
   \begin{equation}
   \Gamma_{\mathfrak{c}}(\mathfrak{n}):=\left\{
   \left(\begin{array}{cc} a & b\\ c & d \end{array} \right)\in GL_2(F)| a,d\in\mathcal{O},b\in \mathfrak{c}^{-1}\mathfrak{d}^{-1}, c\in\mathfrak{ncd}, ad-bc\in U_F^+\right\},
   \end{equation}
   we have
   \begin{equation*}
     (ad-bc)^{k/2}(cz+d)^{-k}f(\frac{az+b}{cz+d})=\varphi(a)f(z).
   \end{equation*}

    Here we are using the same shorthand as in ~\cite[\S 1]{Sh}.  The modularity condition implies that $f$ has a Fourier expansion
    \begin{equation*}
f(z) = a(0) + \sum_{b\in \mathfrak{c}^+} a(b)q^b
    \end{equation*}
where $q^b=e^{2\pi i Tr_{F/\Q}(bz)}$.

 The space of such forms is finite dimensional; we denote this space by $M_{k,\mathfrak{c},\psi}(\mathfrak{n},\C)$. More generally, for any ring $R\subset \C$, let $M_{k,\mathfrak{c},\varphi}(\mathfrak{n},R)$ denote the subset of forms with Fourier coefficients in $R$. Shimura has shown that $M_{k,\mathfrak{c},\varphi}(\mathfrak{n},\overline{\Q})\otimes \C = M_{k,\mathfrak{c},\varphi}(\mathfrak{n},\C)$; using the embedding $\overline{\Q}\hookrightarrow \overline{\Qp}$ fixed at the beginning, we can define $M_{k,\mathfrak{c},\varphi}(\mathfrak{n},R)$ for any subring of $\overline{\Qp}$.

We define a Hilbert modular form (without the $\mathfrak{c}$) of weight $k$ and level $\mathfrak{n}$ to be a $|Cl^+(F)|$-tuple of $\mathfrak{c}$-Hilbert modular forms, where $\mathfrak{c}$ ranges over a set of representatives of the strict ideal classes. We will usually write this as $f = (f_\mathfrak{c})_\mathfrak{c}$. For a ray class character $\chi$ of conductor dividing $\mathfrak{n}$, we will say $f$ has character $\chi$ if $S(\mathfrak{a})f = \chi(\mathfrak{a})$ for almost all prime ideals $\mathfrak{a}$ of $F$ (see ~\cite[p. 648]{Sh}). We denote the space of such forms by $M_{k}(\mathfrak{n},\chi)$.

Given a Hilbert modular form $f$ in the latter sense, the normalized Fourier expansion is defined as follows: for each nonzero integral ideal $\mathfrak{m}$, there is a unique $\mathfrak{c}$ in our choice of strict ideal class representatives for which we can write $\mathfrak{m}\mathfrak{c}=(b)$ for some totally positive $b\in \mathfrak{c}$. Then we let
\begin{equation*}
  c(\mathfrak{m},f) = a_\mathfrak{c}(b)N\mathfrak{c}^{-k/2}.
\end{equation*}
where $a_\mathfrak{c}(b)$ is the coefficient of $q^b$ in $f_\mathfrak{c}$. For each $\lambda\in Cl^+(F)$, we also set
\begin{equation*}
  c_\lambda(0,f) = a_\mathfrak{c}(0) N\mathfrak{c}^{-k/2}.
\end{equation*}
As the notation suggests, neither of these expressions depend on our choice of $b$.

For each prime $\ell\nmid \mathfrak{n}$, and each prime $\mathfrak{q}|\mathfrak{n}$, there are Hecke operators $T_\ell$ and $U_\mathfrak{q}$ which act on the spaces $M_{k}(\mathfrak{n},\chi)$. Fix a rational prime $p$ and suppose $\mathfrak{p}|\mathfrak{n}$ for all $\mathfrak{p}|p$. If $R$ is a complete subring of $\overline{\Qp}$, then we say $f$ is ordinary if $ef=f$, where $e:=\lim_n \prod_{\mathfrak{p}|p} U_\mathfrak{p}^{n!}$.

   Let $\mathfrak{m}_\Lambda$ be the maximal ideal of $\Lambda$, $\Lambda_{(0)}$ its localization at the prime ideal $(T)$, and $F_\Lambda$ its field of fractions. For $E$ a finite extension of $\Qp$, let $\Lambda_E=E\otimes_{\Zp}\Lambda$. Fix an integral ideal $\mathfrak{n}$, and an odd ray class character $\chi$ of conductor dividing $\mathfrak{n}$. Following Wiles, we define $\mathscr{M}_\Lambda^{ord}(\mathfrak{n},\chi)$, the space of level $\mathfrak{n}$ ordinary $\Lambda$-adic forms of character $\chi$, to be a collection of coefficients
  \begin{equation*}
  \{c_\lambda(0,\mathscr{F})\},\{c(\mathfrak{m},\mathscr{F})\}\in \Lambda
  \end{equation*}
  where $\lambda$ runs over $Cl^+(F)$ and $\mathfrak{m}$ runs over the nonzero integral ideals of $\mathcal{O}_F$, such that for almost all pairs $k\geq 2, \zeta \in \mu_{p^\infty}$, the reduction of this system modulo the ideal $P_{\zeta,k}=(T+1-\zeta u^k)$ gives the normalized Fourier coefficients of an ordinary parallel weight $k$ Hilbert modular form of level $\overline{\mathfrak{n}}:=\text{lcm}(p\cdot ord(\zeta),\mathfrak{n})$ and character $\psi_\zeta\chi\omega^{1-k}$. We call the reduction mod $P_k:=P_{1,k}$ the weight $k$ specialization. For any subalgebra $\Lambda \subset R \subset F_\Lambda$, we define $\mathscr{M}_R^{ord}(\mathfrak{n},\chi)=\mathscr{M}_\Lambda^{ord}(\mathfrak{n},\chi)\otimes_\Lambda R$. We also let $\mathscr{S}^{ord}_\Lambda(\mathfrak{n},\chi),\mathscr{S}^{ord}_{F_\Lambda}(\mathfrak{n},\chi),$ etc. denote the corresponding spaces of cusp forms.

  Let $\mathbb{T}^{ord}$ denote the ordinary $\Lambda$-adic Hecke algebra of level $\mathfrak{n}$ and character $\chi$, i.e. the $\Lambda$-algebra generated, for $\ell\nmid \mathfrak{n}p$ and $\mathfrak{q}|\mathfrak{n}p$, by the Hecke operators $T_\ell, U_\mathfrak{q}$ acting on $\mathscr{M}^{ord}_{\Lambda}(\mathfrak{n},\chi)$. Formulae for this action in terms of $q$-expansions are given at the top of page 537 in~\cite{Wi:88}. The following lemma is probably well known, but as far as we know has not been written down.

  \begin{lma} If the weight $k$ specializations of a collection $\{c_\lambda(0,\mathscr{F})\},\{c(\mathfrak{m},\mathscr{F})\}\in \Lambda$ give a classical ordinary form for infinitely many $k\geq 2$, then they are classical for all but finitely many $k\geq 2$.
  \begin{proof}
  Fix a weight $k$ and let $E$ be a finite extension of $\Qp$ containing all the Hecke eigenvalues appearing in $M_k^{ord}(\mathfrak{n}_S,\Qp,\chi\omega^{1-k})$. Wiles shows in ~\cite[Thm. 1.4.1]{Wi:88} that the system of Hecke eigenvalues corresponding to any classical eigenform of weight $k$ can be realized as a quotient of the $\Lambda$-adic Hecke algebra $\mathbb{T}^{ord}/(1+T-u^k) \rightarrow E$. If $\mathfrak{m}$ is the corresponding maximal ideal of $\mathbb{T}^{ord}\otimes E$, the space $\mathscr{M}^{ord}_{\Lambda_{E}}(\mathfrak{n},\chi)/\mathfrak{m}\mathscr{M}^{ord}_{\Lambda_{E}}(\mathfrak{n},\chi)$ is nonzero by Nakayama's lemma. Therefore, any classical eigenform can be realized as the weight $k$ specialization of some $\Lambda_E$-adic form (i.e. we don't need to take a finite extension of $\Lambda_E$; note that we are \emph{not} claiming there is an eigenform defined over $\Lambda_E$).

  We claim that this implies that any ordinary form with coefficients in $\Qp$ can be realized as the specialization of some $\Lambda_{\Qp}$-adic form. Indeed, on the weight $k$ fiber, we can write the ordinary form as an $E$-linear combination of eigenforms. Lifting this linear combination to $\Lambda_E$ gives a form specializing to the one we need, but a priori only has coefficients in $\Lambda_E$. However, by averaging over Gal$(\Lambda_E/\Lambda_{\Qp})$, and observing that specialization intertwines this action with the action of Gal$(E/\Qp)$, we get a $\Lambda_{\Qp}$-adic form with the desired specialization.

The claim implies that the map
\begin{equation*}
\mathscr{M}_{\Lambda_{\Qp}}^{ord}(\mathfrak{n},\chi)/(1+T-u^k)\mathscr{M}_{\Lambda_{\Qp}}^{ord}(\mathfrak{n},\chi)\rightarrow M_k(\mathfrak{n}_S,\Qp, \chi\omega^{1-k})
 \end{equation*}
 is an isomorphism for almost all $k$, say $k\geq k_0$. In particular,
 \begin{equation*}
 \text{rank}_{\Lambda_{\Qp}} \mathscr{M}^{ord}_{\Lambda_{\Qp}}(\mathfrak{n},\chi)= \text{dim}_{\Qp}(\overline{\mathfrak{n}},\Qp,\chi\omega^{1-k})
  \end{equation*}
  for almost all $k$; call this dimension $d$. We can choose ideals $\mathfrak{a}_1,\ldots,\mathfrak{a}_d$ of $\mathcal{O}_F$ so that the map
\begin{align*}
\pi: \mathscr{M}_{\Lambda_{\Qp}}^{ord}&(\mathfrak{n},\chi)\rightarrow (\Lambda_{\Qp})^d\\
\mathscr{F}&\mapsto \left(c(\mathfrak{a}_i,\mathscr{F})\right)_i
 \end{align*}
 is injective. After inverting a finite set of primes $\mathfrak{S}$ of $\Lambda_{\Qp}$, $\pi$ is an isomorphism. Therefore, for $P_k\not\in \mathfrak{S}$, we have
 \begin{align*}
\pi_k: M_k^{ord}&(\overline{\mathfrak{n}},\Qp,\chi\omega^{1-k})\cong \Qp^d\\
f&\mapsto \left(c(\mathfrak{a}_i,f)\right)_i
 \end{align*}

 Now suppose we had a collection of coefficients $\{c_\lambda(0,\mathscr{H})\},\{c(\mathfrak{m},\mathscr{H})\}$ with infinitely many classical specializations. There is a unique element $\mathscr{F}$ of  $\mathscr{M}_{\Lambda_{\Qp}[\frac{1}{\mathfrak{S}}]}^{ord}(\mathfrak{n},\chi)$ such that $c(\mathfrak{a}_i,\mathscr{F})=c(\mathfrak{a}_i,\mathscr{H})$ for all $i$. Moreover, at each weight $k$ with $P_k\notin\mathfrak{S}, k>k_0,$ and where $\mathscr{H}$ is classical, the reduction of $\mathscr{H}$ must agree with the reduction of $\mathscr{F}$ by the isomorphism $\pi_k$. Thus, $\mathscr{H}$ and $\mathscr{F}$ must be equal since they agree on a Zariski dense set. This proves the lemma.

  \end{proof}
  \end{lma}

  A typical example of ordinary $\Lambda$-adic forms of tame level $\mathfrak{n}$ and character $\chi$ are the $\Lambda$-adic Eisenstein series $\mathscr{E}(\eta,\psi)$ attached to a pair of (not necessarily primitive) narrow ray class characters $\eta,\psi$ such that $\eta\psi=\chi$, $\text{cond}(\eta)\text{cond}(\psi)=p\mathfrak{n}$, $(p,\text{cond}(\eta))=1$:
  \begin{equation*}
    c_\lambda(0,\mathscr{E}(\eta,\psi))=\delta_{\eta}2^{-g}\eta^{-1}(\mathfrak{c}_\lambda)\mathcal{L}_{\{\mathfrak{n}p\infty\},\eta^{-1}\psi\omega},
  \end{equation*}
  \begin{equation*}
    c(\mathfrak{m},\mathscr{E}(\eta,\psi))=\sum_{\substack{\mathfrak{r}|\mathfrak{m}\\(\mathfrak{r},p)=1}} \eta(\frac{\mathfrak{m}}{\mathfrak{r}})\psi(\mathfrak{r})\langle N\mathfrak{r}^{-1}\rangle(1+T)^\frac{\log\langle N\mathfrak{r}\rangle}{\log u}.
  \end{equation*}
  Here $\delta_\eta=1$ if cond$(\eta)=1$, and is zero otherwise.

  We denote $G_\zeta:=\mathcal{L}_{\{\mathfrak{p}|p\},1}$, so that $G_\zeta(u^s-1)=\zeta_{F,p}(1-s)$. Let $\mathscr{G}:=2^gG_\zeta^{-1}\mathscr{E}(1,\omega^{-1})$, so that the constant term of $\mathscr{G}$ at each infinite cusp is identically one. It follows from a result of Colmez~\cite{Co:88} that if Leopoldt's Conjecture is true for $(F,p)$ then $G_\zeta$ has a pole of order one at $T=0$. In this case, the form $\mathscr{G}$, which \emph{a priori} only lies in $\mathscr{M}^{ord}_{F_\Lambda}(1,\omega^{-1})$, actually lies in $\mathscr{M}^{ord}_{\Lambda_{(0)}}(1,\omega^{-1})$, with specialization equal to the constant form \textbf{1}, i.e. $c_\lambda(0,\mathscr{G}(0))=1$ and $c(\mathfrak{m},\mathscr{G}(0))=0$ for all $\lambda,\mathfrak{m}$. Theorem 2 in Section 4 below shows that even if Leopoldt fails, there is a suitable cusp form $\mathscr{J}$ such that $\mathscr{G}-\mathscr{J}$ has all of these properties.\\
  \\

\subsection{Conjecture 1.ii}

We now recall the proof of Conjecture \emph{1.ii} given in [DDP], making a few of our own cosmetic changes, but also assuming the existence of the form $\mathscr{J}$ in order to remove Leopoldt from their hypotheses. To prepare for section 3, we will assume that there is a unique prime $\mathfrak{p}$ above $p$ in $F$, as this will highlight how the $\mathscr{L}$-invariant hypothesis comes into play. When there is more than one prime above $p$, the arguments we give showing that one can replace $\mathscr{G}$ by $\mathscr{G}-\mathscr{J}$ in their proof go through unchanged.

Let $F,p,\chi,S,R$ be as in Conjecture 1. The first step in their proof is to obtain a Galois theoretic interpretation of $\mathscr{L}_{alg}$.

Let $E$ be an extension of $\Qp$ containing the values of all characters of conductor dividing $\text{cond}(\chi)p\infty$, and $E(\chi^{-1})$ the $E[Gal(\overline{F}/F)]$-module which is one-dimensional over $E$ and on which Galois acts by $\chi^{-1}$. For ease of notation, let $\Lambda$ denote $\Lambda_E$ (so $p$ is invertible in $\Lambda$). Finally, let $H^1_\mathfrak{p}(F,E(\chi^{-1}))$ be the subspace of $H^1(F,E(\chi^{-1}))$ which is unramified at all primes away from $\mathfrak{p}$, and at $\mathfrak{p}$, lies in the $E$-linear span of $\kappa_{ur}$ and $\kappa_{cyc}$, where $\kappa_{ur}\in H^1(F_\mathfrak{p},E(\chi^{-1})) = Hom(\text{Gal}(\overline{F_\mathfrak{p}}/F_\mathfrak{p}),E)$ is the unramified (additive!) character $\text{Frob}_\mathfrak{p}\mapsto 1$, and $\kappa_{cyc}$ is the restriction of the global character
\begin{equation*}
\text{Gal}(\overline{F}/F)\twoheadrightarrow\text{Gal}(F_\infty/F) \hookrightarrow 1+p\mathbb{Z}_p \stackrel{\log_p}{\rightarrow} \Zp \hookrightarrow E.
\end{equation*}
In the sequel, we will also use $\kappa_{cyc}$ to denote the global character.

In ~\cite[\S 1]{DDP:09}, it is shown that $\dim_EH^1_\mathfrak{p}(F,E(\chi^{-1}))=1$, and that the unique class (up to a scalar) is ramified at $\mathfrak{p}$. In other words, if we write its restriction to $p$ as $x\kappa_{ur} + y\kappa_{cyc}$, then $y\neq 0$. In fact we have
\begin{prp}[\emph{loc. cit.}, Prop. 1.6]
\begin{equation*}
\frac{x}{y} = -\mathscr{L}_{alg}(\chi).
\end{equation*}
\end{prp}
The idea now is to use modular forms to explicitly construct a class in $H^1_\mathfrak{p}(F,E(\chi^{-1}))$ whose restriction to $G_\mathfrak{p}$ can be shown to be equal (up to a scalar) to $-\mathscr{L}_{an}(\chi)\kappa_{ur} + \kappa_{cyc}$.

Denote by $\chi_R$ the character of conductor $R$ which has the same primitive as $\chi$. Consider the level $R$ weight one Hilbert modular Eisenstein series $E_1(1,\chi_R)$. We have
\begin{equation*}
c_\lambda(0,E_1(1,\chi_R))=2^{-g}L_R(\chi,0)+\delta_\chi \chi^{-1}(\lambda)L_R(\chi^{-1},0),
\end{equation*}

where $\delta_\chi=1$ if cond$(\chi)=1$, and is 0 otherwise. We also have
\begin{equation*}
 U_\mathfrak{p}E_1(1,\chi_R)=E_1(1,\chi_R)+E_1(1,\chi_S),
  \end{equation*}
  which implies that $e=\lim_n(U_\mathfrak{p}^{n!})$ acts by the identity on this form. Thus, the ordinary $\Lambda$-adic form $\mathscr{P}^0:=e[(\mathscr{G}-\mathscr{J})E_1(1,\chi_R)]$ has weight one specialization equal to $E_1(1,\chi_R)$. Moreover, its constant terms satisfy
  \begin{equation*}
  c_\lambda(0,\mathscr{P}^0)=2^{-g}L_R(\chi,0)+\delta_\chi \chi^{-1}(\lambda)L_R(\chi^{-1},0)
  \end{equation*}
   independent of the weight, since this is clearly true before taking the ordinary projection, and the only Eisenstein series contributing to these cusps at the classical higher weight specializations are already ordinary.

Over $F_\Lambda$, we can decompose $\mathscr{P}^0$ into a linear combination of a cusp form and ordinary Eisenstein series. The coefficients $a(1,\chi)$ and $a(\chi,1)$ of the Eisenstein families $\mathscr{E}(1,\chi)$ and $\mathscr{E}(\chi,1)$ in this decomposition are computed in \emph{loc. cit.} $\S 2$, using knowledge of the constant terms of $e[\mathscr{G}E_1(1,\chi_R)]$ at \emph{all} unramified cusps, not just the infinite cusps. The weight $k$ specializations of these coefficients are given respectively by
\begin{align*}
a(1,\chi)(k)&=\frac{L_R(\chi,0)}{L_{S,p}(\chi\omega,1-k)}=\frac{-1}{\mathscr{L}_{an}(\chi,k)}\\
a(\chi,1)(k)&=\frac{L_R(\chi^{-1},0)\langle N\mathfrak{n}\rangle^{k-1}}{L_{S,p}(\chi^{-1}\omega,1-k)}=\frac{-\langle N\mathfrak{n}\rangle^{k-1}}{\mathscr{L}_{an}(\chi^{-1},k)}.
\end{align*}
Note that these computations are unaffected if one replaces $\mathscr{G}$ by $\mathscr{G}-\mathscr{J}$ since $\mathscr{J}$ vanishes at all cusps.

The Eisenstein series other than $\mathscr{E}(1,\chi)$ and $\mathscr{E}(\chi,1)$ can be killed by an appropriate application of Hecke operators away from $p$ without affecting the weight 1 specialization $E_1(1,\chi_R)$, simply by dividing each Hecke operator by its eigenvalue on $E_1(1,\chi)$. It follows that there is some $t$, a linear combination of Hecke operators away from $p$, acting by the identity on $E_1(1,\chi)$, such that
\begin{equation*}
t(\mathscr{P}^0-a(1,\chi)\mathscr{E}(1,\chi)-a(\chi,1)\mathscr{E}(\chi,1))
\end{equation*}
 is a cusp form. Note that $a(1,\chi)$ and $a(\chi,1)$ have poles at weight one of order equal to the order of vanishing of the corresponding $p$-adic $L$-functions; let us say they are of order $d_\chi$ and $d_{\chi^{-1}}$, respectively. The $\mathscr{L}$-invariant condition is equivalent to
\begin{equation*}
 ord_{k=1}a(1,\chi)=ord_{k=1}(a(1,\chi)+a(\chi,1)).
  \end{equation*}

  Let $\pi\in \Lambda$ be the uniformizer at weight one given by $\frac{1}{u\log u}(1+T-u)$. Then the universal cyclotomic character $\boldsymbol{\chi}_{cyc}$ can be written $\boldsymbol{\chi}_{cyc}=1+\kappa_{cyc}\pi+O(\pi^2)$.

 Since the poles of $a(1,\chi)$ and $a(\chi,1)$ do not cancel, we have
 \begin{equation*}
   \mathscr{F}=\mathscr{P}^0-a(1,\chi)\mathscr{E}(1,\chi)-a(\chi,1)\mathscr{E}(\chi,1)\in \pi^{-d_\chi}\mathscr{M}^{ord}_{\Lambda}(\mathfrak{n},\chi).
 \end{equation*}

 We want to consider the image of this form in
  \begin{equation*}
  \pi^{-d_\chi}\mathscr{M}^{ord}_{\Lambda}(\mathfrak{n},\chi)/\pi^{-d\chi+2}\mathscr{M}^{ord}_{\Lambda}(\mathfrak{n},\chi)
  \end{equation*}
and compute the Hecke action. To do this, we will use the following identities
 \begin{align*}
  \mathscr{E}(\chi,1)\equiv \mathscr{E}(1,\chi) &\equiv E_1(1,\chi_S) (\text{mod }\pi)\\
  \mathscr{P}^0 &\equiv E_1(1,\chi_R)(\text{mod }\pi)\vspace{0.2cm}\\
  T_\ell \mathscr{E}(\chi,1) = (\chi(\ell)+\boldsymbol{\chi}_{cyc}(\ell))\mathscr{E}(\chi,1)\quad&\quad U_\mathfrak{p} \mathscr{E}(\chi,1)= \mathscr{E}(\chi,1)\\
  T_\ell \mathscr{E}(1,\chi) = (1+\chi(\ell)\boldsymbol{\chi}_{cyc}(\ell))\mathscr{E}(1,\chi) \quad&\quad U_\mathfrak{p} \mathscr{E}(1,\chi) = \mathscr{E}(1,\chi)\\
  T_\ell E_1(1,\chi_R) = (1+\chi(\ell))E_1(1,\chi_R) \quad&\quad U_\mathfrak{p} E_1(1,\chi_R) = E_1(1,\chi_R) + E_1(1,\chi_S)
  \end{align*}
  We compute:

   \begin{align*}
     T_\ell\mathscr{F} &= (1+\chi(\ell))\mathscr{P}^0 - (1+\chi(\ell)\boldsymbol{\chi}_{cyc}(\ell))a(1,\chi)\mathscr{E}(1,\chi) - (\boldsymbol{\chi}_{cyc}(\ell) + \chi(\ell))a(\chi,1)\mathscr{E}(\chi,1))\\
     &=(1+\chi(\ell))\mathscr{F} - \chi(\ell)\kappa_{cyc}(\ell)\pi a(1,\chi)\mathscr{E}(1,\chi) - \kappa_{cyc}(\ell)\pi a(\chi,1)\mathscr{E}(\chi,1)\\
     &= (1+\chi(\ell))\mathscr{F} - (\chi(\ell)\kappa_{cyc}(\ell)\pi a(1,\chi) + \kappa_{cyc}(\ell)\pi a(\chi,1))E_1(1,\chi_S)\\
     &= (1+\chi(\ell))\mathscr{F} + (\chi(\ell)\kappa_{cyc}(\ell)a(1,\chi) + \kappa_{cyc}(\ell)a(\chi,1)) (\frac{\pi}{a(1,\chi)+a(\chi,1)})\mathscr{F}\\
     &= \left(\left(1+ \pi \frac{a(\chi,1)}{a(1,\chi)+a(\chi,1)}\kappa_{cyc}(\ell)\right) + \chi(\ell)\left(1+ \pi \frac{a(1,\chi)}{a(1,\chi)+a(\chi,1)}\kappa_{cyc}(\ell)\right)\right)\mathscr{F};
   \end{align*}

\begin{align*}
  U_\mathfrak{p}\mathscr{F} &= \mathscr{F} + E_1(1,\chi_S)= (1- \frac{1}{a(1,\chi)+a(\chi,1)})\mathscr{F}=(1+\pi\frac{a(1,\chi)}{a(1,\chi)+a(\chi,1)}\mathscr{L}_{an})\mathscr{F}.
\end{align*}

The calculation for the $U_\mathfrak{q}$ operators is similar.

Thus we get a map
\begin{equation*}
\mathbb{T}^{ord}\rightarrow \Lambda/\pi^2 \cong E[\pi]/\pi^2
\end{equation*}
\begin{align*}T_\ell &\mapsto (1+\pi\frac{a(\chi,1)}{a(1,\chi)+a(\chi,1)}\kappa_{cyc}(\ell))+ \chi(\ell)(1+\pi\frac{a(1,\chi)}{a(1,\chi)+a(\chi,1)}\kappa_{cyc}(\ell))\\
U_\mathfrak{q}&\mapsto 1-\pi\frac{a(\chi,1)}{a(1,\chi)+a(\chi,1)}\kappa_{cyc}(\mathfrak{q})\\
U_\mathfrak{p}&\mapsto 1+\pi\frac{a(1,\chi)}{a(1,\chi)+a(\chi,1)}\mathscr{L}_{an}.
\end{align*}

Let $I$ denote the kernel of this map, and $\mathfrak{m}$ the maximal ideal containing $I$. Since the image of $t$ is a unit, this map factors through the cuspidal quotient of $\mathbb{T}^{ord}$. Let $R$ denote the localization at $\mathfrak{m}$ of this cuspidal quotient, and think of $I$ and $\mathfrak{m}$ as being ideals of $R$. Finally, let $F_R:=R\otimes_\Lambda F_\Lambda$.

There is a Galois representation
\begin{align*}
\rho:G_F\rightarrow GL_2(F_R)\\
\rho(\sigma)=\left(\begin{array}{cc} a_\sigma & b_\sigma\\ c_\sigma & d_\sigma \end{array} \right).
\end{align*}
unramified at all $\ell\nmid p\mathfrak{n}$ such that Trace$(Frob_\ell)=T_\ell$. For some choice of complex conjugation $c\in G_F$, we may assume $\rho(c)=\left(\begin{array}{cc} 1 & 0 \\ 0 & -1 \end{array} \right)$.

By a standard argument, the $R$-submodule of $F_R$ generated by all $a_\sigma$ is $R$, and similarly for all $d_\sigma$. Moreover,
\begin{align*}
a_\sigma\equiv 1+\pi\frac{a(\chi,1)}{a(1,\chi)+a(\chi,1)}\kappa_{cyc}(\sigma)(\text{mod }I)\\
d_\sigma\equiv \chi(\sigma)(1+\pi\frac{a(1,\chi)}{a(1,\chi)+a(\chi,1)}\kappa_{cyc}(\sigma))(\text{mod }I).
\end{align*}
 By a theorem of Wiles, there is a change-of-basis matrix $\left(\begin{array}{cc} A_\mathfrak{p} & B_\mathfrak{p}\\ C_\mathfrak{p} & D_\mathfrak{p} \end{array} \right)$ with the property that
\begin{equation*}
\left(\begin{array}{cc} a_\sigma & b_\sigma\\ c_\sigma & d_\sigma \end{array} \right)\left(\begin{array}{cc} A_\mathfrak{p} & B_\mathfrak{p}\\ C_\mathfrak{p} & D_\mathfrak{p} \end{array} \right) =
\left(\begin{array}{cc} A_\mathfrak{p} & B_\mathfrak{p}\\ C_\mathfrak{p} & D_\mathfrak{p} \end{array} \right) \left(\begin{array}{cc} \boldsymbol{\chi}_{cyc}\eta_\mathfrak{p}^{-1}(\sigma) & *\\ 0 & \eta_\mathfrak{p}(\sigma) \end{array} \right)
\end{equation*}
for all $\sigma\in G_\mathfrak{p}$. Here $\eta_\mathfrak{p}$ is the unramified character $\text{Frob}_\mathfrak{p}^k\mapsto U_\mathfrak{p}^k$.

We note that in the basis with complex conjugation diagonalized, $\rho(G_\mathfrak{p})$ is \emph{not} upper-triangular on any component of $F_R$.\footnote{If it were, the function $c_\sigma$ reduced mod $\mathfrak{m}$ would yield a non-trivial element of $H^1(G_F,E(\chi))$ which is unramified everywhere (recall that $\mathfrak{p}$ is the unique prime above $p$). However, there are no unramified elements in this $H^1$, since e.g. there are no unramified $\Zp$-extensions of $F_\chi$.} Therefore, $C_\mathfrak{p}$ is invertible in $F_R$. Hence, for $\sigma\in G_\mathfrak{p}$, we can write \begin{equation*}
b_\sigma=\frac{A_\mathfrak{p}}{C_\mathfrak{p}}[\boldsymbol{\chi}_{cyc}\eta_\mathfrak{p}^{-1}(\sigma) - a_\sigma].
\end{equation*}

If $B$ is the $R$-module generated by $b_\sigma$ (or equivalently by $\frac{b_\sigma}{d_\sigma}$) as $\sigma$ ranges over $G_F$, then $B$ is finite over $R$ by a standard compactness argument, so that $B/\mathfrak{m}B$ is nonzero. If $K:G_F\rightarrow B/\mathfrak{m}B$ is the composition of $b/d$ with reduction mod $\mathfrak{m}$, then $[K]\in H^1(G_F,B/\mathfrak{m}B(\chi^{-1}))$ is unramified outside $\mathfrak{p}$ and nonzero, since if $K$ were a coboundary, one can check it would have to be identically zero by considering $K(c)$. Since there are no everywhere unramified elements of this $H^1$, $[K]$ must be ramified at $\mathfrak{p}$; in particular, it is nontrivial at $\mathfrak{p}$. This argument, combined with Nakayama's lemma, shows that $B$ is in fact generated by $b_\sigma$ for $\sigma\in G_\mathfrak{p}$.

Reducing the above equation modulo $(I\frac{A_\mathfrak{p}}{C_\mathfrak{p}}\cap B)$, we get for $\sigma\in G_\mathfrak{p}$,

\begin{align*}
\overline{b(\sigma)}&=\frac{A_\mathfrak{p}}{C_\mathfrak{p}}[\boldsymbol{\chi}_{cyc}(\sigma)\eta_\mathfrak{p}^{-1}(\sigma)-a(\sigma)]\\
&=\frac{A_\mathfrak{p}}{C_\mathfrak{p}}[\boldsymbol{\chi}_{cyc}(\sigma)U_\mathfrak{p}^{-\kappa_{ur}(\sigma)}-a(\sigma)]\\
&=\frac{A_\mathfrak{p}}{C_\mathfrak{p}}[(1+\pi\kappa_{cyc}(\sigma))(1-\pi\frac{a(1,\chi)}{a(1,\chi)+a(\chi,1)}\mathscr{L}_{an}\kappa_{ur}(\sigma))-1-\pi\frac{a(\chi,1)}{a(1,\chi)+a(\chi,1)}\kappa_{cyc}(\sigma)]\\
&=\frac{A_\mathfrak{p}}{C_\mathfrak{p}}\pi[\frac{a(1,\chi)}{a(1,\chi)+a(\chi,1)}\kappa_{cyc}(\sigma)-\frac{a(1,\chi)}{a(1,\chi)+a(\chi,1)}\mathscr{L}_{an}\kappa_{ur}(\sigma)].
\end{align*}

Since $\pi\in\mathfrak{m}$ and $\mathfrak{m}^2\subset I$, the module $B/(I\frac{A_\mathfrak{p}}{C_\mathfrak{p}}\cap B)$ is $\mathfrak{m}$-torsion. Furthermore, by our initial assumptions we can choose $\sigma$ so that the bracketed expression is a unit. It follows that $B/(I\frac{A_\mathfrak{p}}{C_\mathfrak{p}}\cap B)$ is cyclic over $R/\mathfrak{m}\cong E$, generated by $\pi\frac{A_\mathfrak{p}}{C_\mathfrak{p}}$. Since it is $\pi$-torsion, we must have $\overline{b(\sigma)}=\overline{b/d(\sigma)}$. Thus, we may view $\overline{b/d}$ as a cocycle with coefficients in $E(\chi^{-1})$, whose restriction to $G_\mathfrak{p}$ is given, after dividing by the unit $\frac{a(1,\chi)}{a(1,\chi)+a(\chi,1)}$, by $\kappa_{cyc}-\mathscr{L}_{an}\kappa_{ur}$. By Proposition 1, $\mathscr{L}_{alg}=\mathscr{L}_{an}$, and the proof is complete.

\section{Removing the $\mathscr{L}$-invariant condition}
Let $\mathscr{L}_\chi$ denote the reciprocal of the leading term of $a(1,\chi)$, and $\mathscr{L}_{\chi^{-1}}$ the reciprocal of the leading term of $a(\chi,1)$, with respect to the uniformizer $\pi$. Then $\mathscr{L}_\chi=-\mathscr{L}_{an}(\chi)$ \emph{only if} $d_\chi=1$; otherwise $\mathscr{L}_{an}(\chi)=0$ (similarly for $\chi^{-1}$).

As explained in the introduction, we break into two cases.

\subsection{$d_\chi < d_{\chi^{-1}}$}

In this case, the form $\mathscr{P}^0 - a(1,\chi)\mathscr{E}(1,\chi) - a(\chi,1)\mathscr{E}(\chi,1)$ has a pole of order $d_{\chi^{-1}}$ coming from the third term, while the second term has a pole of order $d_\chi$. Consider the Hecke action on the image of this form in $\pi^{-d_{\chi^{-1}}}\mathscr{M}^{ord}_{\Lambda}(\mathfrak{n},\chi) /\pi^{-d\chi+2}\mathscr{M}^{ord}_{\Lambda}(\mathfrak{n},\chi)$. A similar analysis to the previous case shows that this form is an eigenform, giving rise to a homomorphism
\begin{equation*}
  \mathbb{T}^{ord}\rightarrow E[\pi]/\pi^{d_{\chi^{-1}}-d_\chi+2},
\end{equation*}
as follows:
\begin{align*}T_\ell &\mapsto \left(\boldsymbol{\chi}_{cyc}(\ell)-\pi\frac{a(1,\chi)}{a(1,\chi)+a(\chi,1)}\kappa_{cyc}(\ell)\right) + \chi(\ell)\left(1+\pi\frac{a(1,\chi)}{a(1,\chi)+a(\chi,1)}\kappa_{cyc}(\ell)\right) \\
U_\mathfrak{q}&\mapsto 1-\pi\frac{a(\chi,1)}{a(1,\chi)+a(\chi,1)}\kappa_{cyc}(\mathfrak{q})\\
U_\mathfrak{p}&\mapsto 1 - \frac{1}{a(\chi,1)}.
\end{align*}

As before, this homomorphism factors through $R$, the localization of the cuspidal quotient at the maximal ideal corresponding to the $E_1(1,\chi)$-system of Hecke eigenvalues. Let $I\subset R$ denote the kernel. There is a Galois representation $G_F\rightarrow GL_2(F_R)$ with coefficients $a_\sigma,b_\sigma,c_\sigma,d_\sigma$, satisfying
\begin{align*}
R\ni a_\sigma& \equiv \boldsymbol{\chi}_{cyc}(\sigma)-\pi\frac{a(1,\chi)}{a(1,\chi)+a(\chi,1)}\kappa_{cyc}(\sigma)(\text{mod }I)\\
R\ni d_\sigma& \equiv \chi(\sigma)(1+\pi\frac{a(1,\chi)}{a(1,\chi)+a(\chi,1)}\kappa_{cyc}(\sigma))(\text{mod }I).
\end{align*}
Finally, there is a transition matrix $\left(\begin{array}{cc} A_\mathfrak{p} & B_\mathfrak{p}\\ C_\mathfrak{p} & D_\mathfrak{p} \end{array} \right)$ which satisfies
\begin{equation*}
b_\sigma = \frac{A_\mathfrak{p}}{C_\mathfrak{p}}[\boldsymbol{\chi}_{cyc}\eta_\mathfrak{p}^{-1}(\sigma)-a(\sigma)]
\end{equation*}
for all $\sigma\in G_\mathfrak{p}$. Reducing this equation modulo $I\frac{A_\mathfrak{p}}{C_\mathfrak{p}}\cap B$, we get for $\sigma\in G_\mathfrak{p}$,
\begin{align*}
\overline{b/d(\sigma)}=\overline{b(\sigma)}=& \frac{A_\mathfrak{p}}{C_\mathfrak{p}}[\boldsymbol{\chi}_{cyc}\eta_\mathfrak{p}^{-1}(\sigma)-a(\sigma)]\\
=& \frac{A_\mathfrak{p}}{C_\mathfrak{p}}[\frac{1}{a(\chi,1)}\kappa_{ur}(\sigma) + \pi\frac{a(1,\chi)}{a(\chi,1)}\kappa_{cyc}(\sigma)]\\
=& \frac{A_\mathfrak{p}}{C_\mathfrak{p}}[\pi^{d_{\chi^{-1}}}\mathscr{L}_{\chi^{-1}}\kappa_{ur}(\sigma)+ \pi^{d_{\chi^{-1}}-d_\chi+1}\frac{\mathscr{L}_{\chi^{-1}}}{\mathscr{L}_{\chi}}\kappa_{cyc}(\sigma)].
\end{align*}

For $\sigma\in I_\mathfrak{p}\backslash \ker(\kappa_{cyc})$, the bracketed expression generates $(\pi^{d_{\chi^{-1}}-d_\chi+1}\frac{A_\mathfrak{p}}{C_\mathfrak{p}})$ mod $I\frac{A_\mathfrak{p}}{C_\mathfrak{p}}$. Hence, as before, $B/(B\cap I\frac{A_\mathfrak{p}}{C_\mathfrak{p}})$ is one-dimensional over $R/\mathfrak{m}\cong E$, with a canonical generator given by $\pi^{d_{\chi^{-1}}-d_\chi+1}\frac{A_\mathfrak{p}}{C_\mathfrak{p}}$. The composition $\overline{b/d}:G_F\rightarrow B/(B\cap I\frac{A_\mathfrak{p}}{C_\mathfrak{p}})\cong E$ given by this generator yields a nonzero cocycle
\begin{equation*}
[\kappa]\in H^1_\mathfrak{p}(F,E(\chi^{-1}))
\end{equation*}
with the property that
\begin{equation*}
[\kappa]|_{G_\mathfrak{p}}=\delta_{d_\chi=1}\mathscr{L}_{\chi^{-1}}\kappa_{ur} +\frac{\mathscr{L}_{\chi^{-1}}}{\mathscr{L}_{\chi}}\kappa_{cyc}.
\end{equation*}
After multiplying by the nonzero scalar $\frac{\mathscr{L}_{\chi}}{\mathscr{L}_{\chi^{-1}}}$, the right hand side equals $-\mathscr{L}_{an}(\chi)\kappa_{ur} + \kappa_{cyc}$. This finishes the proof.
\\
\\
\subsection{$d_\chi = d_{\chi^{-1}}$ and $\mathscr{L}_\chi=-\mathscr{L}_{\chi^{-1}}$}

In this case, although the second and third term of $\mathscr{P}^0 - a(1,\chi)\mathscr{E}(1,\chi) - a(\chi,1)\mathscr{E}(\chi,1)$ each have poles of order $d=d_\chi$, the sum only has a pole of order $d-1$. We consider the Hecke action on the image of this form in $\pi^{-d+1}\mathscr{M}^{ord}_{\Lambda}(\mathfrak{n},\chi) /\pi^{-d+2}\mathscr{M}^{ord}_{\Lambda}(\mathfrak{n},\chi)$. Unlike the previous cases, this form is \emph{not} an eigenform. Nevertheless, we have
\begin{align*}
&T_\ell (\mathscr{P}^0 - a(1,\chi)\mathscr{E}(1,\chi) - a(\chi,1)\mathscr{E}(\chi,1))=\\ &(1+\chi(\ell))\mathscr{P}^0-(1+\chi(\ell)(1+\pi\kappa_{cyc}(\ell)))a(1,\chi)\mathscr{E}(1,\chi) - (1+\pi\kappa_{cyc}(\ell)+\chi(\ell))a(\chi,1)\mathscr{E}(\chi,1)=\\
&(1+\chi(\ell))[\mathscr{P}^0 - a(1,\chi)\mathscr{E}(1,\chi) - a(\chi,1)\mathscr{E}(\chi,1)] - (\kappa_{cyc}(\ell)\chi(\ell)\pi a(1,\chi)+ \kappa_{cyc}(\ell)\pi a(\chi,1))E_1(1,\chi).\\
\\
&U_\mathfrak{q}(\mathscr{P}^0 - a(1,\chi)\mathscr{E}(1,\chi) - a(\chi,1)\mathscr{E}(\chi,1))=\\
&(\mathscr{P}^0 - a(1,\chi)\mathscr{E}(1,\chi) - a(\chi,1)\mathscr{E}(\chi,1)) -\pi\kappa_{cyc}(\mathfrak{q})a(\chi,1)E_1(1,\chi).\\
\\
&U_\mathfrak{p}(\mathscr{P}^0 - a(1,\chi)\mathscr{E}(1,\chi) - a(\chi,1)\mathscr{E}(\chi,1))=\\
&(\mathscr{P}^0 - a(1,\chi)\mathscr{E}(1,\chi)- a(\chi,1)\mathscr{E}(\chi,1)) + E_1(1,\chi).
\end{align*}

Thus, although the image of our form is not an eigenvector for the Hecke operators, it \emph{is} a generalized eigenvector for the $E_1(1,\chi)$-system of eigenvalues; the Hecke stable subspace it generates is two-dimensional over $E$, with a basis given by $(\mathscr{P}^0 - a(1,\chi)\mathscr{E}(1,\chi) - a(\chi,1)\mathscr{E}(\chi,1))$ and $\pi^{-d+1}E_1(1,\chi)$. The Hecke action can then be viewed as a homomorphism
\begin{equation*}
\mathbb{T}^{ord}\rightarrow \left(\begin{array}{cc} x & y\\ 0 & x \end{array} \right)\subset M_2(E).
\end{equation*}
The image is canonically isomorphic to $E[\varepsilon]/\varepsilon^2$. Under this identification, the map $\mathbb{T}^{ord}\rightarrow E[\varepsilon]/\varepsilon^2$ is given explicitly by
\begin{align*}
T_\ell&\mapsto 1+\chi(\ell) - \varepsilon(\frac{\kappa_{cyc}(\ell)\chi(\ell)}{\mathscr{L}_\chi}+ \frac{\kappa_{cyc}(\ell)}{\mathscr{L}_{\chi^{-1}}})\\
U_\mathfrak{q}&\mapsto 1 - \varepsilon(\frac{\kappa_{cyc}(\mathfrak{q})}{\mathscr{L}_{\chi^{-1}}})\\
U_\mathfrak{p}&\mapsto 1 +\varepsilon(\delta_{d=1}).
\end{align*}
Following the same proof as before, we get a Galois representation $G_F\rightarrow GL_2(F_R)$ such that
\begin{align*}
R\ni a_\sigma &\equiv 1 - \varepsilon\frac{\kappa_{cyc}(\sigma)}{\mathscr{L}_{\chi^{-1}}}(\text{mod }I)\\
R\ni d_\sigma &\equiv \chi(\sigma)(1-\varepsilon\frac{\kappa_{cyc}(\ell)}{\mathscr{L}_\chi})(\text{mod }I).
\end{align*}
However, in this case the image of the universal cyclotomic character $\boldsymbol{\chi}_{cyc}$ in $R/I$ is trivial, as we are working ``purely in weight one.'' Thus, reducing the equation
\begin{equation*}
b_\sigma = \frac{A_\mathfrak{p}}{C_\mathfrak{p}}[\boldsymbol{\chi}_{cyc}\eta_\mathfrak{p}^{-1}-a_\sigma]
\end{equation*}
modulo $I\frac{A_\mathfrak{p}}{C_\mathfrak{p}}$ gives for $\sigma\in G_\mathfrak{p}$
\begin{align*}
\overline{b/d(\sigma)}&= \frac{A_\mathfrak{p}}{C_\mathfrak{p}}[\eta_\mathfrak{p}^{-1}-a_\sigma]\\
&= \frac{A_\mathfrak{p}}{C_\mathfrak{p}}\varepsilon(-\delta_{d=1}\kappa_{ur}(\sigma)+\frac{\kappa_{cyc}(\sigma)}{\mathscr{L}_{\chi^{-1}}}).
\end{align*}
Just as before, we see that $B/(B\cap I\frac{A_\mathfrak{p}}{C_\mathfrak{p}})$ is one-dimensional, generated by $\varepsilon\frac{A_\mathfrak{p}}{C_\mathfrak{p}}$, and that the function $\overline{b/d}$ yields a class $[\kappa]\in H^1_\mathfrak{p}(F,E(\chi^{-1}))$ such that
\begin{align*}
[\kappa]|_{G_\mathfrak{p}} &= -\delta_{d=1}\kappa_{ur} + \frac{1}{\mathscr{L}_{\chi^{-1}}}\kappa_{cyc}\\
&= -\delta_{d=1}\kappa_{ur} - \frac{1}{\mathscr{L}_{\chi}}\kappa_{cyc}.
\end{align*}
Multiplying through by $-\mathscr{L}_\chi$, the right hand side becomes
\begin{equation*}
\delta_{d=1}\mathscr{L}_\chi\kappa_{ur}+\kappa_{cyc}=-\mathscr{L}_{an}(\chi)\kappa_{ur}+\kappa_{cyc}.
 \end{equation*}
 This finishes the proof.

\section{A $\Lambda$-adic form passing through 1}

  We revert to letting $\Lambda$ denote $\Zp[[T]]$. In this section we prove the following
  \begin{thm}
  There exists an $F_\Lambda$-adic cusp form $\mathscr{J}\in \mathscr{S}^{ord}_{F_\Lambda}(1,\omega^{-1})$ such that $\mathscr{G}-\mathscr{J}\in \mathscr{M}^{ord}_{\Lambda_{(0)}}(1,\omega^{-1})$ and $(\mathscr{G}-\mathscr{J})(0)$ is the constant form \textbf{1}.
  \end{thm}

The existence of this form is perhaps well known to experts. As explained in Section 2, this theorem removes the reliance on Leopoldt's conjecture in the proof of Conjecture $1.ii$. At the end of this section, we also explain how the form $\mathscr{J}$ gives an easier and more direct construction of the Iwasawa extensions corresponding to the ``Leopoldt'' zeros in Wiles' proof of the Main Conjecture.
\subsection{Reductions and geometric $\Lambda$-adic forms}
We begin with some reductions. First, it is enough to find any ordinary $\Lambda_{E,(0)}$-adic form of level one with constant weight zero specialization, where $E$ is a finite extension of $\Qp$. Rescaling, we may assume the constant at weight zero is equal to one. Let $\mathscr{F}$ be such a form, so that
\begin{equation*}
\mathscr{F}\in \mathscr{M}^{ord}_{F_{\Lambda_E}}(1)=\bigoplus_{\chi}\mathscr{M}^{ord}_{F_{\Lambda_E}}(1,\chi\omega^{-1})
\end{equation*}
where $\chi$ ranges over all even ray class characters of conductor $1$. It causes no harm to assume $E$ contains sufficiently many roots of unity, so that there is a $E$-linear combination of diamond operators projecting onto the $\omega^{-1}$ component above. Specializing at weight 0, this acts by the identity on the constant form, so we may assume $\mathscr{F}$ has character $\omega^{-1}$. We can write $\mathscr{F}$ as an $F_{\Lambda_W}$-linear combination of a cusp form and Eisenstein series $\mathscr{E}(\eta,\eta^{-1}\omega^{-1})$, where $\eta$ ranges over strict ray class characters of conductor 1. For $\eta\neq 1$, the system of Hecke eigenvalues associated to the Eisenstein series $\mathscr{E}(\eta,\eta^{-1}\omega^{-1})$ at weight zero differs from that of the constant form. Hence, there is a $\Lambda_{E,(0)}$-linear combination of Hecke operators which will kill all Eisenstein contributions except $\mathscr{E}(1,\omega^{-1})$, and act by the identity on the constant form. Applying this to $\mathscr{F}$, we are left with an $F_{\Lambda_E}$-linear combination of $\mathscr{E}(1,\omega^{-1})$ and a cusp form. Since the constant terms of this form are identically one, it must be equal to $\mathscr{G}-\mathscr{J}$ for some $F_{\Lambda_E}$-cuspform $\mathscr{J}$. Finally, we can average over Gal$(F_{\Lambda_E}/F_\Lambda)$ so that $\mathscr{J}$ has coefficients in $F_\Lambda$.

Our next reduction requires a new definition. Let $r>0$ be a natural number, and $\mathbf{\Lambda}:=\Zp[[\frac{T}{p^r}]]$. For an integral ideal $\mathfrak{n}$ and odd ray class character $\chi$ of conductor dividing $\mathfrak{n}$, we define an ordinary $\mathbf{\Lambda}$-adic Hilbert modular form to be a collection of coefficients
\begin{equation*}
  \{c_\lambda(0,\mathscr{F})\},\{c(\mathfrak{m},\mathscr{F})\}\in \mathbf{\Lambda}
  \end{equation*}
  such that for infinitely many $k\in p^r\mathbb{N}$, their image under the specialization $T=u^k-1$ are the coefficients of a classical ordinary Hilbert modular form of parallel weight $k$, level lcm$(p,\mathfrak{n})$, and character $\chi\omega^{1-k}$. We will denote this module by $\mathscr{M}_{\mathbf{\Lambda}}^{ord}(\mathfrak{n},\chi)$, for any ring $\mathbf{\Lambda}\subset R\subset F_\mathbf{\Lambda}$, we set $\mathscr{M}_{R}^{ord}(\mathfrak{n},\chi)= \mathscr{M}_{\mathbf{\Lambda}}^{ord}(\mathfrak{n},\chi)\otimes_{\mathbf{\Lambda}} R$. Using the constant dimensionality of the spaces of ordinary weight $k$ forms, and an argument similar to that used in Lemma 2, one can show that the space of $F_\mathbf{\Lambda}$-adic forms has the same dimension as the space of $F_\Lambda$-adic forms, and hence is identified with $\mathscr{M}^{ord}_{F_\Lambda} \otimes_{F_\Lambda} F_\mathbf{\Lambda}$. Now suppose we can find an $\mathbf{\Lambda}_{(\frac{T}{p^r})}$-adic form specializing to $\mathbf{1}$ at weight zero (i.e. modulo $\frac{T}{p^r}$). Then writing it as an $F_\mathbf{\Lambda}$-linear combination of elements of $\mathscr{M}^{ord}_\Lambda$, we can replace the coefficients in $F_\mathbf{\Lambda}$ with elements of $F_\Lambda$ having the same principal part and constant term at weight zero, to arrive at an $F_\Lambda$-adic form specializing to $\mathbf{1}$ at weight zero. Thus, it is enough to find an ordinary family with coefficients in $F_\mathbf{\Lambda}$.

To construct a family of level one, we first construct a family of some auxiliary level $\mathfrak{q}$ using powers of a certain theta series (Lemmas 3 and 4). We will then use the Atkin-Lehner operators $U_\mathfrak{q}$ and $W_\mathfrak{q}$ to project the form down to level one. In order to define these operators, we will make use of a geometric description of $\Lambda$-adic forms (see Proposition 2).

  Fix $F, p$ as in Section 2 and let $\mathfrak{n}$ be an integral ideal of $\mathcal{O}_F$. Fix a strict ideal class, and let $\mathfrak{c}$ be a prime-to-$p$ representative of this class. Let $R$ be a $p$-adically complete DVR. Following ~\cite[Definition 3.2]{AG:05}, we let $T_{m,n}=\mathfrak{M}(R/p^m,\mu_{p^n},\Gamma_0(\mathfrak{n}))$ be the moduli stack over $R/p^m$ whose objects over $S$, for any $R/p^n$-scheme $S$, are given by isomorphism classes of tuples $(A,\iota,\lambda,\phi_\mathfrak{n},i_{p^n})$ where

\begin{itemize}
\item $A\rightarrow S$ is an abelian scheme of relative dimension $g$;
\item $\iota: \mathcal{O}_F\hookrightarrow End_R(A)$ is a ring homomorphism;
\item $\lambda: (M_A,M_A^+) \stackrel{\cong}{\rightarrow}(\mathfrak{c},\mathfrak{c}^+)$
is an $\mathcal{O}_F$-linear isomorphism of $\acute{\text{e}}$tale sheaves over $T$ between the module of symmetric $\mathcal{O}_F$-linear homomorphisms from $A$ to its dual $A^\vee$ to the ideal $\mathfrak{c}$, such that the polarizations $M_A^+$ map to $\mathfrak{c}^+$;
\item $\phi_{\mathfrak{n}}\subseteq A$ is an $\mathcal{O}_F$-invariant closed subgroup scheme which is isomorphic to the constant group scheme $(\mathcal{O}_F/\mathfrak{n})$ $\acute{\text{e}}$tale locally on $S$;
\item $i_{p^n}: \mu_{p^n}\otimes_{\Z}\mathfrak{d}^{-1}\hookrightarrow A$ is an inclusion of group schemes.
\end{itemize}

  These are referred to as $\mathfrak{c}$-polarized Hilbert Blumenthal Abelian Varieties (HBAV's) with level structure. Following ~\cite[Definition 11.4]{AG:05}, we define a $p$-adic $\mathfrak{c}$-Hilbert modular form (or $\mathfrak{c}$-HMF) of level $\Gamma_0(\mathfrak{n})$ over $R$ to be an element of
\begin{equation*}
  V_{\infty,\infty}:=\displaystyle\underleftarrow{\lim}_m\underrightarrow{\lim}_n H^0(T_{m,n/(R/p^m)},\mathcal{O}_{T_{m,n}}).
\end{equation*}

If $\chi: (\mathcal{O}_F\otimes \Zp)^\times\rightarrow R^\times$ is a finite order character, we will say the form is of (parallel) weight $k\in \Zp$ and character $\chi$ if for any $\alpha\in (\mathcal{O}_F\otimes \Zp)^\times$, we have
\begin{equation*}
\alpha^*(f) = \chi(\alpha)\mathbf{Nm}(\alpha)^k f,
\end{equation*}
where $\alpha^*f(A,\iota,\lambda,\phi_\mathfrak{n},i_{p^\infty})=f(A,\iota,\lambda,\phi_\mathfrak{n},i_{p^\infty}\circ\alpha^{-1})$, and $\textbf{Nm}:(\mathcal{O}_F\otimes \Zp)^\times\rightarrow 1+2p\Zp\rightarrow R^\times$ is induced by the norm map followed by projection onto the $1$-units.

Fix an isomorphism $\epsilon: \mathfrak{c}\otimes \Zp\cong \mathcal{O}_F\otimes \Zp$. The \emph{$q$-expansion at $\infty$} of $f$ is an element
\begin{equation*}
  f(q)\in R[[q^b]]_{b\in \mathfrak{c}^+\cup \{0\}}.
\end{equation*}
which generalizes the $q$-expansion of classical Hilbert modular forms.  We refer to ~\cite[Definition 11.6]{AG:05} for the precise definition. In their notation, it is the evaluation of $f$ at the cusp $(\mathfrak{c},\mathcal{O}_F,\epsilon,j_\epsilon)$ where $j_\epsilon$ is induced from $\epsilon$ as in [\emph{loc. cit.}, 6.5]. The \emph{q-expansion principle} states that a $p$-adic $\mathfrak{c}$-HMF of weight $\kappa$ is determined by its $q$-expansion [\emph{loc. cit.}, 11.7].

  For every $k\geq 2$, there is a Hecke-equivariant inclusion $M_{k,\mathfrak{c}}(\Gamma_0(\mathfrak{n}p),R)\hookrightarrow V_{\infty,\infty/R}$ which preserves $q$-expansions and weights~\cite[Thm. 1.10.15]{Ka:78}. In fact, in the quoted theorem, the space of classical forms on the left hand side of the inclusion is more general than the forms we considered in Section 2: it allows any power of $p$ in the level, and in the complex setting, it consists of those forms invariant under the subgroup of $\Gamma_{\mathfrak{c}}(\mathfrak{n})$ consisting of matrices of determinant 1 (see~\cite[6.11]{AG:05}). However, this certainly contains the forms we want to consider, and this is all we will need. We will call a form in the image of the above inclusion classical.

Let $W$ be a finite flat DVR over $\Zp$. Let $\mathfrak{m}_\Lambda$ denote the maximal ideal of $\Lambda_W$. We now present two definitions of ``$p$-adic'' $\Lambda_W$-adic forms and prove that they are the same. We also want similar statements to hold for $\mathbf{\Lambda}$-adic forms, and will indicate where changes need to be made.

 Recall that $u\in1+2p\Zp$ is a generator of the image of $\text{Gal}(F_\infty/F)$. Define the map $\phi:(\mathcal{O}_F\otimes \Zp)^\times \rightarrow \Lambda_W^\times$ by the composition
\begin{equation*}
  (\mathcal{O}_F\otimes \Zp)^\times\rightarrow \text{Gal}(F_\infty/F) \stackrel{u\mapsto 1+T}{\longrightarrow}  \Lambda_W^\times
\end{equation*}
We may also consider $\phi$ as a $\mathbf{\Lambda}_W$-valued character via the inclusion $\Lambda\subset \mathbf{\Lambda}$.

\textbf{Definition.} A \emph{Wiles} $\Lambda_W$-adic $\mathfrak{c}$-Hilbert modular form $\mathscr{F}$, of level $\Gamma_0(\mathfrak{n})$, is a multiset of elements $\{c(t,\mathscr{F})\}_{t\in\mathfrak{c}^+\cup \{0\}}\subset \Lambda_W$ such that for every $s\in \Zp$, the sequence of elements of $W$ obtained from the specialization $T\mapsto u^s-1$ is the $q$-expansion of a $p$-adic $\mathfrak{c}$-Hilbert modular form of level $\Gamma_0(\mathfrak{n})$ and weight $s$ over $W$. We use the same definition for $\mathbf{\Lambda}_W$-adic $\mathfrak{c}$-Hilbert modular forms, except that we only require the specialization condition to hold for $s\in p^r\Zp$.

\indent \emph{Remark.} Given a $|Cl^+(F)|$-tuple of Wiles $\Lambda_W$-adic forms, one for each strict ideal class with representative $\mathfrak{c}$, and with infinitely many specializations giving the Fourier coefficients of an ordinary classical Hilbert modular form, we can obtain a $\Lambda$-adic form as in Section 2 as follows. Under the usual normalization for weight $k$ forms, one would set $c(\mathfrak{m},\mathscr{F}(u^k-1))=(N\mathfrak{c})^{-k/2}\cdot a_\mathfrak{c}(b)(u^k-1)$, where $b$ is a totally positive generator of $\mathfrak{m}\mathfrak{c}$. However, this presents a problem since $(N\mathfrak{c})^{-k/2}$ may not vary $p$-adically continuously with $k$. So instead, we simply set $c(\mathfrak{m},\mathscr{F})= a_\mathfrak{c}(b)$, and $c_\mathfrak{c}(0,\underline{\quad})$ equal to the constant term. Since this is independent of the choice of $b$ at infinitely many weights, it must be independent $\Lambda$-adically. This modification will not affect what we are ultimately interested in: finding a family whose weight zero specialization is the constant form $\textbf{1}$.

\textbf{Definition.} A \emph{Katz} $\Lambda_W$-adic $\mathfrak{c}$-Hilbert modular form of level $\Gamma_0(\mathfrak{n})$ is an element of the subspace of $\displaystyle V_{\infty,\infty/\Zp}\widehat{\otimes}_{\Zp}\Lambda_W:=\underleftarrow{\lim}V_{\infty,\infty/\Zp}\otimes_{\Zp}\Lambda_W/\mathfrak{m}_\Lambda^n$ satisfying
\begin{equation}\tag{*}
\alpha^*(f)= \phi(\alpha)f.
 \end{equation}
for every $\alpha\in (\mathcal{O}_F\otimes \Zp)^\times$. We use the same definition for $\mathbf{\Lambda}_W$ mutatis mutandis.

This last property is equivalent to requiring that for every $s\in \Zp$ (resp. $p^r\Zp$), reducing the form modulo $(1+T-u^s)$ yields a $p$-adic $\mathfrak{c}$-Hilbert modular form of weight $s$ defined over $\Lambda_W/(1+T-u^s)$. Note that neither of these definitions require any of the specializations to be classical.

By definition, a Katz $\Lambda$-adic $\mathfrak{c}$-HMF is nothing but a compatible sequence of $p$-adic $\mathfrak{c}$-HMFs over $\Lambda_W/\mathfrak{m}_\Lambda^n$ satisfying certain extra conditions. Thus, we may define the $q$-expansion at $\infty$ of a Katz $\Lambda$-adic form to be the inverse limit of these $q$-expansions; it is an element of $\Lambda[[q^b]]_{b\in \mathfrak{c}^+\cup\{0\}}$.

The following proposition is due to Hida when $F=\Q$~\cite[Thm. 3.2.16]{Hi:00}. We essentially follow his proof.

\begin{prp} The space of Katz $\Lambda_W$-adic $\mathfrak{c}$-Hilbert modular forms is identified with the space of Wiles $\Lambda_W$-adic $\mathfrak{c}$-Hilbert modular forms via $q$-expansion at $\infty$. The same is true for $\mathbf{\Lambda}_W$-adic forms.
\begin{proof}
We first explain the proof for $\Lambda_W$-adic forms. It follows from the definitions that the $q$-expansion of a Katz form is a Wiles $\Lambda_W$-adic form, so we need only show that all Wiles $\Lambda_W$-adic forms arise in this way. Let $\mathscr{F}$ be a Wiles $\Lambda_W$-adic form. We start by reinterpreting $\mathscr{F}$ as a measure $C(\Zp,\Zp)\rightarrow V_{\infty,\infty/W}$, defined by sending the function ${x}\choose{n}$ to the coefficient of $T^n$ in $\mathscr{F}$ (which is a $p$-adic $\mathfrak{c}$-HMF by virtue of being a limit of $p$-adic $\mathfrak{c}$-HMF's). By the binomial theorem, this measure has the property that for $s\in \Zp$, the function $x\mapsto u^{sx}$ is sent to $\mathscr{F}(u^s-1)$, a $p$-adic $\mathfrak{c}$-HMF of weight $s$.

Taking the completed tensor product with $\Lambda$ of this measure gives a map
\begin{equation*}
C(\Zp,\Lambda)\rightarrow V_{\infty,\infty/W}\hat{\otimes}\Lambda.
\end{equation*}
The image of the function $x\mapsto (1+T)^x$ is easily seen to be a Katz $\Lambda_W$-adic form (i.e. obeys the equation (*)), with $q$-expansion equal to $\mathscr{F}$.

Now suppose $\mathscr{F}$ is a Wiles $\mathbf{\Lambda}_W$-adic form. Define the submodule $M\subset C(\Zp,\Zp)$ by demanding that if $C(\Zp,\Zp)\ni f = \sum_{N\geq 0} a_N \binom{x}{N}$, then
\begin{equation*}
  f\in M \iff \frac{a_N}{p^{rN}}\in \Zp \text{ and } \frac{a_N}{p^{rN}}\rightarrow 0\text{ as }N\rightarrow \infty.
\end{equation*}
Then we may consider $\mathscr{F}$ as a measure $M\rightarrow V_{\infty,\infty/W}$ which sends $p^{rN}\binom{x}{N}$ to the coefficient of $(\frac{T}{p^r})^N$ in $\mathscr{F}$. Then just as before, we take the completed tensor product with $\mathbf{\Lambda}_W$, and consider the image of the function $x\mapsto (1+T)^x$. This gives the desired Katz $\mathbf{\Lambda}_W$-adic form.

\end{proof}
\end{prp}

The Katz definition gives us a geometric interpretation of $\Lambda$-adic forms as follows:

Let $\mathbb{M}_\mathfrak{c}$ denote the functor from the category of $\mathfrak{m}_\Lambda$-adically complete $\Lambda$-algebras to Sets which takes an algebra $R$ to the set of isomorphism classes tuples $(A,\iota,\lambda,\phi_\mathfrak{n},i_{p^\infty})$ as above. Then we can view a $\Lambda$-adic $\mathfrak{c}$-HMF as a natural transformation from this functor to the forgetful functor $\mathbb{A}^1$, which further satisfies $(*)$. For $\mathbf{\Lambda}$-adic $\mathfrak{c}$-HMFs, the same statement holds if we consider $\mathfrak{m}_\mathbf{\Lambda}$-adically complete $\mathbf{\Lambda}$-algebras.\\
 \\
 \subsection{Construction of the form} We now construct the level one ordinary family with constant weight zero specialization. We begin by quoting Lemma 1.4.2 of~\cite{Wi:88}, which is attributed to Hida:

\begin{lma}
For some prime $\mathfrak{q}\nmid p$, and some $m>0$, there is a Hilbert modular form $f$ of weight $2^m(p-1)$ and level $\Gamma_0(p\mathfrak{q})$, with coefficients in $\Zp$, such that $c_\lambda(0,f)=1$ for all $\lambda$, and $f\equiv 1(\text{mod }p)$.
\end{lma}

\emph{Remark}. The quoted lemma has a power $p^j$ in the level; however, we can apply the operator $U_p:=T_{\mathfrak{n}p}(p)$ $j-1$ times to decrease the level at $p$ to $\Gamma_0(p)$ without altering any of the other properties. This lemma is proved using theta series coming from the extension $F(\mu_p)/F(\mu_p)^+$. An alternative approach is to use lifts of suitable powers of the Hasse invariant (see~\cite[Lemma 11.10]{AG:05}).

Write $f=(f_\mathfrak{c})_\mathfrak{c}$. Since $f_\mathfrak{c}\equiv 1(\text{mod }p)$, for any $s\in \Zp$, we can make sense of $f^s$ as a $p$-adic $\mathfrak{c}$-Hilbert modular form of level $\mathfrak{q}$ and weight $2^m(p-1)s$ (here we are using the equivalence of Katz-type and Serre-type $p$-adic Hilbert modular forms in parallel weight; see ~\cite[Theorem 11.12]{AG:05}).

Let $e$ be the $p$-adic valuation of $u-1$. Let

\begin{displaymath}
   r = \left\{
     \begin{array}{lr}
     e+m+1 &\text{if }p=2\\
     e+1   &\text{if }p>2
     \end{array}
   \right.
\end{displaymath}
and $\mathbf{\Lambda}=\Zp[[\frac{T}{p^r}]]$.

\begin{lma}
There is a $\mathbf{\Lambda}$-adic $\mathfrak{c}$-HMF $\mathscr{F}_\mathfrak{c}$ such that $\mathscr{F}_\mathfrak{c}(u^s-1)=f_\mathfrak{c}^{s/2^m(p-1)}$ for all $s\in 2^m\Zp$.
\begin{proof}
Write $f_\mathfrak{c}=\sum_{b\in\mathfrak{c}^+\cup \{0\}} c_b q^b$. Fix a positive integer $k$, and some $b\in \mathfrak{c}^+$. Let $\Pi$ be the set of all tuples $\{(k_1,b_1),\ldots,(k_\ell,b_\ell)\}$, where $k_i\in \mathbb{N}_{>0}$ and $b_i\in \mathfrak{c}^+$, such that $\sum k_ib_i=b$.
Note that the cardinality of $\Pi$ is finite and does not depend on $k$. The $q^b$ coefficient of $f_\mathfrak{c}^k$ is given explicitly by
\begin{equation*}
[q^b]f_\mathfrak{c}^k = \sum_{\Pi} \binom{k}{k_1,\ldots,k_\ell,k-(k_1+\ldots+k_\ell)}c_{b_{1}}^{k_1}\ldots c_{b_{\ell}}^{k_\ell},
\end{equation*}
where we have used that $c_0=1$, and we interpret a multinomial coefficient with negative arguments to be zero.

As this sum is finite, it suffices to prove that each element of the above sum is given by evaluating some element of $\mathbf{\Lambda}$ at $u^{(p-1)2^mk}-1$. The multinomial coefficient can be written as $\frac{P(k)}{k_1!\ldots k_\ell!}$ for some polynomial $P$. Since $f_\mathfrak{c}\equiv 1(\text{mod }p)$, we have $v_p(c_{b_i})\geq 1$, so $v_p(\frac{c_{b_i}^{k_i}}{k_i!})\geq 0$. Thus, it is enough to show $P(k)$ can be expressed as an element of $\mathbf{\Lambda}$, and for this it is enough to show $k$ itself can be. The function $k$ is nothing but the weight divided by $2^m(p-1)$, so it is given by
\begin{equation*}
\frac{\log_p(1+T)}{(p-1)2^m\log_p u} \in \mathbf{\Lambda}.
\end{equation*}
This concludes the proof.
\end{proof}
\end{lma}

Note that the weight zero specialization of the tuple $\mathscr{F}=(\mathscr{F}_\mathfrak{c})_\mathfrak{c}$ is the constant form \textbf{1}, and that this form has infinitely many classical specializations. To remove $\mathfrak{q}$ from the level, we make use of the Hecke operators $W_\mathfrak{q}$ and $U_\mathfrak{q}$, interpreting $\mathscr{F}$ as a rule on $\mathfrak{c}$-polarized HBAV's for some $\mathfrak{c}$.

 For $C$ a finite subgroup scheme of an HBAV $A$, we let $\pi_C:A\rightarrow A/C$. We can geometrically define the operators $U_\mathfrak{q}$ and $W_\mathfrak{q}$ on $\Lambda_E$ (or $\mathbf{\Lambda}_E$)-adic modular forms of level $\Gamma_0(\mathfrak{q})$ in the usual way (see e.g. ~\cite[pg. 320-321]{Hi:06}):
\begin{equation*}
  U_\mathfrak{q}\mathscr{F}(A,\lambda,i_{p^{\infty}},\phi_q)=\frac{1}{N\mathfrak{q}}\sum_{C\cap \phi_\mathfrak{q}=\{0\}} \mathscr{F}(A/C,\pi_{C*}\lambda,\pi_C\circ i_{p^\infty},\pi_{C*}\phi_q)
\end{equation*}
\begin{equation*}
  W_\mathfrak{q}\mathscr{F}(A,\lambda,i_{p^{\infty}},\phi_\mathfrak{q})=
  \mathscr{F}(A/\phi_\mathfrak{q},\pi_{\phi_\mathfrak{q}*}\lambda,\pi_{\phi_\mathfrak{q}}\circ i_{p^\infty},\pi_{\phi_\mathfrak{q}*}A[\mathfrak{q}]),
\end{equation*}

The $F_{\mathbf{\Lambda}_{E}}$-adic form $e\left(\frac{W_\mathfrak{q}+U_\mathfrak{q}}{1+N\mathfrak{q}^{-1}}\right)\mathscr{F}$ is ordinary of level one, since its evaluation at any tuple $(A,\lambda,i_{p^{\infty}},\phi_q)$ does not depend on level $\mathfrak{q}$ structure. It has infinitely many classical specializations, and its weight zero specialization is the constant form \textbf{1}. By the remarks at the beginning of the section, this finishes the proof of Theorem 1.
 $\Box$
\subsection{Application to the Iwasawa Main Conjecture}
Theorem 2 allows us to give a direct construction of the Iwasawa extensions corresponding to the (conjecturally nonexistent) zeroes of the $p$-adic zeta function at $s=1$~\cite{Wi:90}\footnote{C. Khare has explained to the author that one can also construct the extra extensions by allowing ramification at an auxiliary prime away from $p$.}. A separate proof is needed for these extensions, as the general argument only constructs a space of extensions of rank $\text{ord}_{s=1}\zeta_{F,p}(s)$, but the Main Conjecture predicts that this space has rank $\delta=\text{ord}_{s=1}\zeta_{F,p}(s)+1$. The proofs given in $\S\S 10,11$ of \emph{loc. cit.} are somewhat indirect, using ``patching'' arguments similar to what is needed in the weight one case. The proof we give here is relatively straightforward with the help of Theorem 2.

Since the non-constant terms of the form \textbf{1} vanish, we have that for each nonzero integral ideal $\mathfrak{m}$,
\begin{equation*} c(\mathfrak{m},\mathscr{G}-\mathscr{J})=c(\mathfrak{m},2^nG_\zeta^{-1}\mathscr{E}(1,\omega^{-1})-\mathscr{J})\in \mathfrak{m}_{(0)}.
\end{equation*}
Thus,
\begin{equation*}
c(\mathfrak{m},\mathscr{E}(1,\omega^{-1})-2^{-n}G_\zeta\mathscr{J})\in \mathfrak{m}_{(0)}^\delta.
\end{equation*}

 Consider the action of the cuspidal Hecke algebra on $2^{-n}G_\zeta\mathscr{J}$. We have
 \begin{equation*}
  2^{-n}G_\zeta\mathscr{J}\equiv\mathscr{E}(1,\omega^{-1})(\text{mod }\mathfrak{m}_{(0)}^\delta)
   \end{equation*}

   \noindent away from the constant terms, so there is a map
 \begin{equation*}
\mathbb{T}^{cusp}\rightarrow \Lambda/\mathfrak{m}_{(0)}^\delta
\end{equation*}
 which is just the $\mathscr{E}(1,\omega^{-1})$ system of Hecke eigenvalues $(\text{mod }\mathfrak{m}_{(0)}^\delta)$. From here, the usual argument by Ribet's method constructs the desired extensions.

\bibliographystyle{alpha}
\bibliography{one}

  \end{document}